\documentclass{gtart_h}

\def\ifplaintex{\expandafter\ifx\csname documentclass\endcsname\relax}

\def\gtp{{\mathsurround=0pt\it $\cal G\mskip-2mu$eometry \&\ 
$\cal T\!\!$opology $\cal P\!$ublications}}  

\def\recd{{\small Received:\qua\receiveddate\ifx\reviseddate\relax
\else\qquad Revised:\qua\reviseddate\fi\par}} 

\def\lognumber#1{\def\thelognumber{#1}}
\def\volumenumber#1{\def\thevolumenumber{#1}}
\def\volumeyear#1{\def\thevolumeyear{#1}}
\def\papernumber#1{\def\thepapernumber{#1}}
\def\pagenumbers#1#2{\def\startpage{#1}\def\finishpage{#2}}
\def\published#1{\def\publishdate{#1}}

\def\received#1{\def\receiveddate{#1}}
\def\revised#1{\def\reviseddate{#1}}
\def\accepted#1{\def\accepteddate{#1}}

\def\asciiaddress#1{\def\theasciiaddress{#1}}

\long\def\asciiabstract#1{\long\def\theasciiabstract{#1}}

\let\\\par\let\thelognumber\relax\let\thevolumenumber\relax
\let\thepapernumber\relax\let\thevolumeyear\relax\let\startpage\relax
\let\finishpage\relax\let\publishdate\relax\let\receiveddate\relax
\let\reviseddate\relax\let\accepteddate\relax\let\theasciititle\relax
\let\theasciiauthors\relax\let\theasciiaddress\relax
\let\theasciiabstract\relax

\let\theasciiemail\relax

\ifplaintex
\font\logobig=cmssbx10 scaled 3836
\font\logomed=cmssbx10 scaled 2557
\else
\font\logobig=cmssbx10 scaled 4200
\font\logomed=cmssbx10 scaled 2800
\fi

\long\def\makeagttitle{   
\count0=\startpage
\agt\hfill      
\hbox to 45truept{\vbox to 0pt{\vglue -13truept{\logomed A\kern -.37em{\logobig 
T}\kern -.38em G}\vss}\hss}
\break
{\small Volume \thevolumenumber\ (\thevolumeyear)
\startpage--\finishpage\nl
Published: \publishdate}

\vglue .25truein

{\parskip=0pt\leftskip 0pt plus
1fil\def\\{\par\smallskip}{\Large\bf\thetitle}\par\medskip} \vglue
0.05truein

{\parskip=0pt\leftskip 0pt plus 1fil\def\\{\par}{\sc\theauthors}
\par\medskip}%
 
\vglue 0.03truein 

{\small\leftskip 25truept\rightskip 25truept{\bf Abstract}\stdspace\theabstract

{\bf AMS Classification}\stdspace\theprimaryclass
\ifx\thesecondaryclass\relax\else; \thesecondaryclass\fi\par
{\bf Keywords}\stdspace \thekeywords\par}\vglue 7truept

}   

\ifplaintex
\font\phead=cmsl9 scaled 950
\font\pnum=cmbx10 scaled 913
\font\pfoot=cmsl9 scaled 950
\headline{\vbox to 0pt{\vskip -4.5mm\line{\small\phead\ifnum
\count0=\startpage ISSN 1472-2739 (on-line) 1472-2747 (printed)
\hfill {\pnum\folio}\else\ifodd\count0\def\\{ }%
\ifx\theshorttitle\relax\thetitle\else\theshorttitle\fi\hfill{\pnum\folio}
\else\def\\{ and }{\pnum\folio}\hfill\ifx\theshortauthors\relax\theauthors
\else\theshortauthors\fi\fi\fi}\vss}}
\footline{\vbox to 0pt{\vglue 0mm\line{\small\pfoot\ifnum\count0=\startpage
\copyright\ \gtp\hfill\else
\agt, Volume \thevolumenumber\ (\thevolumeyear)\hfill\fi}\vss}}
\else
\font=cmsl9 scaled 1050
\font\lnum=cmbx10 
\font=cmsl9 scaled 1050
\makeatletter
\def\@oddhead{{\small\lhead\ifnum\count0=\startpage ISSN 1472-2739 
(on-line) 1472-2747 (printed)\hfill {\lnum\number\count0}\else\ifodd\count0
\def\\{ }\ifx\theshorttitle\relax \thetitle \else\theshorttitle\fi\hfill
{\lnum\number\count0}\else\def\\{ and }{\lnum\number\count0}
\hfill\ifx\theshortauthors\relax 
\theauthors\else\theshortauthors\fi\fi\fi}}\def\@evenhead{\@oddhead}
\def\@oddfoot{\small\lfoot\ifnum\count0=\startpage\copyright\ \gtp\hfill\else
\agt, Volume \thevolumenumber\ (\thevolumeyear)\hfill\fi}
\def\@evenfoot{\@oddfoot}
\makeatother
\fi
\let\maketitlepage\makeagttitle
\let\makeshorttitle\maketitlepage
\let\maketitle\maketitlepage

\newwrite\gtoutfile
\long\gdef\makeheadfile{  
{\def\\{, }\def\s{ }
\immediate\openout\gtoutfile head.xxx
\immediate\write\gtoutfile{Proxy-for: \ifx\theasciiauthors\relax
\theauthors\else\theasciiauthors\fi\s<\ifx\theasciiemail\relax\theemail\else\theasciiemail\fi>}
\immediate\write\gtoutfile{\noexpand\\}
\immediate\write\gtoutfile{Authors: \ifx\theasciiauthors\relax
\theauthors\else\theasciiauthors\fi}
{\def\\{ }\immediate\write\gtoutfile{Title: \ifx\theasciititle\relax
\thetitle\else\theasciititle\fi}}
\immediate\write\gtoutfile{Subj-class: GT or SG, GR etc}
\immediate\write\gtoutfile{MSC-class: \theprimaryclass\ifx\thesecondaryclass\relax\else, \thesecondaryclass\fi}
\immediate\write\gtoutfile{Journal-ref: Algebraic and Geometric Topology \thevolumenumber\s
(\thevolumeyear) \startpage-\finishpage}
\immediate\write\gtoutfile{Comments: Published by Algebraic and
Geometric Topology at}
\immediate\write\gtoutfile{\s\s\s  http://www.maths.warwick.ac.uk/agt/AGTVol\thevolumenumber/agt-\thevolumenumber-\thepapernumber.abs.html}
\immediate\write\gtoutfile{\noexpand\\}
\immediate\write\gtoutfile{}
\ifx\theasciiabstract\relax
\immediate\write\gtoutfile{\theabstract}\else
\immediate\write\gtoutfile{\theasciiabstract}\fi
\immediate\write\gtoutfile{}
\immediate\write\gtoutfile{\noexpand\\}
\immediate\write\gtoutfile{}
\immediate\closeout\gtoutfile}}  

\def\maketitlepage{\makeagttitle\makeheadfile}
\let\makeshorttitle\maketitlepage
\let\maketitle\maketitlepage

\lognumber{9}
\volumenumber{4}
\volumeyear{2004}
\papernumber{9}
\published{23 March 2004}
\pagenumbers{133}{149}
\received{12 August 2003}
\revised{21 January 2004}
\accepted{9 February 2004}

\usepackage{amsmath, amssymb, latexsym}
\usepackage{graphicx, subfigure, epsfig, psfrag}

\def\bkC{{\rm \kern.24em \vrule width.05em height1.4ex depth-.05ex \kern-.26em C}}
\def\C{\bkC}
\def\bksC{{\rm \kern.24em \vrule width.05em height1ex depth-.05ex \kern-.26em C}}

\def\bkE{{\rm I\kern-.22em E}}
    
\def\bkH{{\rm I\kern-.22em H}}
                    
\def\bkQ{{\rm \kern.24em \vrule width.05em height1.4ex depth-.05ex \kern-.26em Q}}

\def\bkR{{\rm I\kern-.17em R}}

\def\bkZ{{\rm Z\kern-.32em Z}}
\def\Z{\bkZ}
\def\bksZ{{\rm Z\kern-.22em Z}}

\def\SL{SL_2(\C)}
\def\PSL{PSL_2(\C)}
\def\F{\mathfrak{F}}
\def\PF{\overline{\mathfrak{F}}}
\def\k{\mathfrak{k}}
\def\X{\mathfrak{X}}
\def\PX{\overline{\mathfrak{X}}}
\def\R{\mathfrak{R}}
\def\PR{\overline{\mathfrak{R}}}
\def\Red{\mathfrak{Red}}
\def\M{\mathfrak{M}}
\def\PM{\overline{\mathfrak{M}}}
\def\T{\X_\tau}
\def\PT{\PX_\tau}
\def\S{\R_\tau}
\def\PS{\PR_\tau}
\def\graph{\mathcal{G}}

\def\tree{\mathcal{T}}

\def\prho{\overline{\rho}}
\def\map{\mu}
\def\pmap{\overline{\mu}}

\def\G{\Gamma}

\DeclareMathOperator{\Hom}{Hom}
\DeclareMathOperator{\tr}{tr}

\DeclareMathOperator{\im}{im}

\DeclareMathOperator{\te}{t}
\DeclareMathOperator{\pte}{\overline{t}}

\DeclareMathOperator{\qe}{q}
\DeclareMathOperator{\identity}{id}

\newtheorem{thm}{Theorem}
\newtheorem{lem}[thm]{Lemma}
\newtheorem{cor}[thm]{Corollary}
\newtheorem{pro}[thm]{Proposition}

\begin{document}
\title{Character varieties of mutative 3--manifolds}
\author{Stephan Tillmann}
\address{D\'epartment de math\'ematiques,
Universit\'e du Qu\'ebec \`a Montr\'eal\\Case 
postale 8888, Succursale Centre-Ville\\Montr\'eal (Qu\'ebec),
Canada H3C 3P8} 
\asciiaddress{Department de mathematiques,
Universite du Quebec a Montreal\\Case postale 8888, 
Succursale Centre-Ville\\Montreal (Quebec), 
Canada H3C 3P8}
\email{tillmann@math.uqam.ca}
\begin{abstract}
We describe a birational map between subvarieties in the character varieties
of mutative 3--manifolds. By studying the birational map, one can decide in
certain circumstances whether a mutation surface is detected by an ideal point
of the character variety. 
\end{abstract}
\asciiabstract{We describe a birational map between subvarieties in the character varieties
of mutative 3-manifolds. By studying the birational map, one can decide in
certain circumstances whether a mutation surface is detected by an ideal point
of the character variety.}

\maketitle


\section{Introduction}

Let $M$ be an irreducible 3--manifold. A \emph{mutation surface} $(S, \tau )$
is a properly embedded incompressible, $\partial$--incompressible surface $S$,
which is not boundary parallel, and an orientation preserving involution
$\tau$ of $S$. The manifold obtained by cutting $M$ open along $S$ and
regluing via $\tau$ is a \emph{$(S,\tau )$--mutant} of $M$, and denoted by
$M^\tau$. Two manifolds are \emph{mutative} if they are related by a finite
sequence of mutations. Mutative manifolds share many geometric and topological
properties. Work of Ruberman \cite{r} and Cooper and Long \cite{a-poly1} shows
that a relationship between the character varieties of mutants exists, and we
now formalise this for both $\SL$-- and $\PSL$--character varieties.

\begin{figure}[t]
\psfrag{a}{{\small $a$}}
\psfrag{b}{{\small $b$}}
\psfrag{c}{{\small $c$}}
\psfrag{d}{{\small $d$}}
\psfrag{abc}{{\small $abc$}}
\psfrag{t}{{\small $\tau$}}
  \begin{center}
    \subfigure[$S_3$]{
      \includegraphics[height=3.2cm]{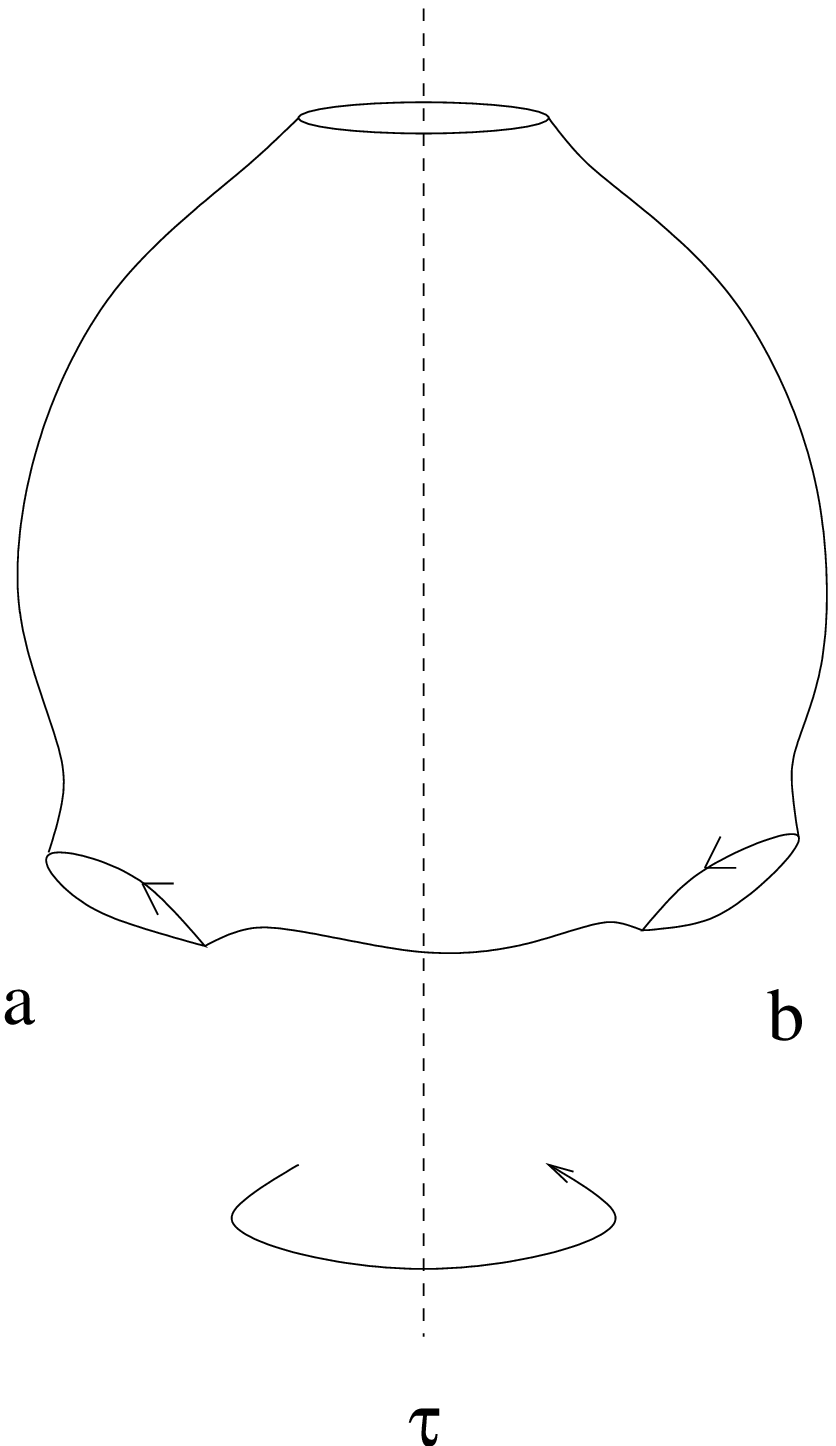}}
      \qquad
    \subfigure[$S_4$]{
      \includegraphics[height=3.2cm]{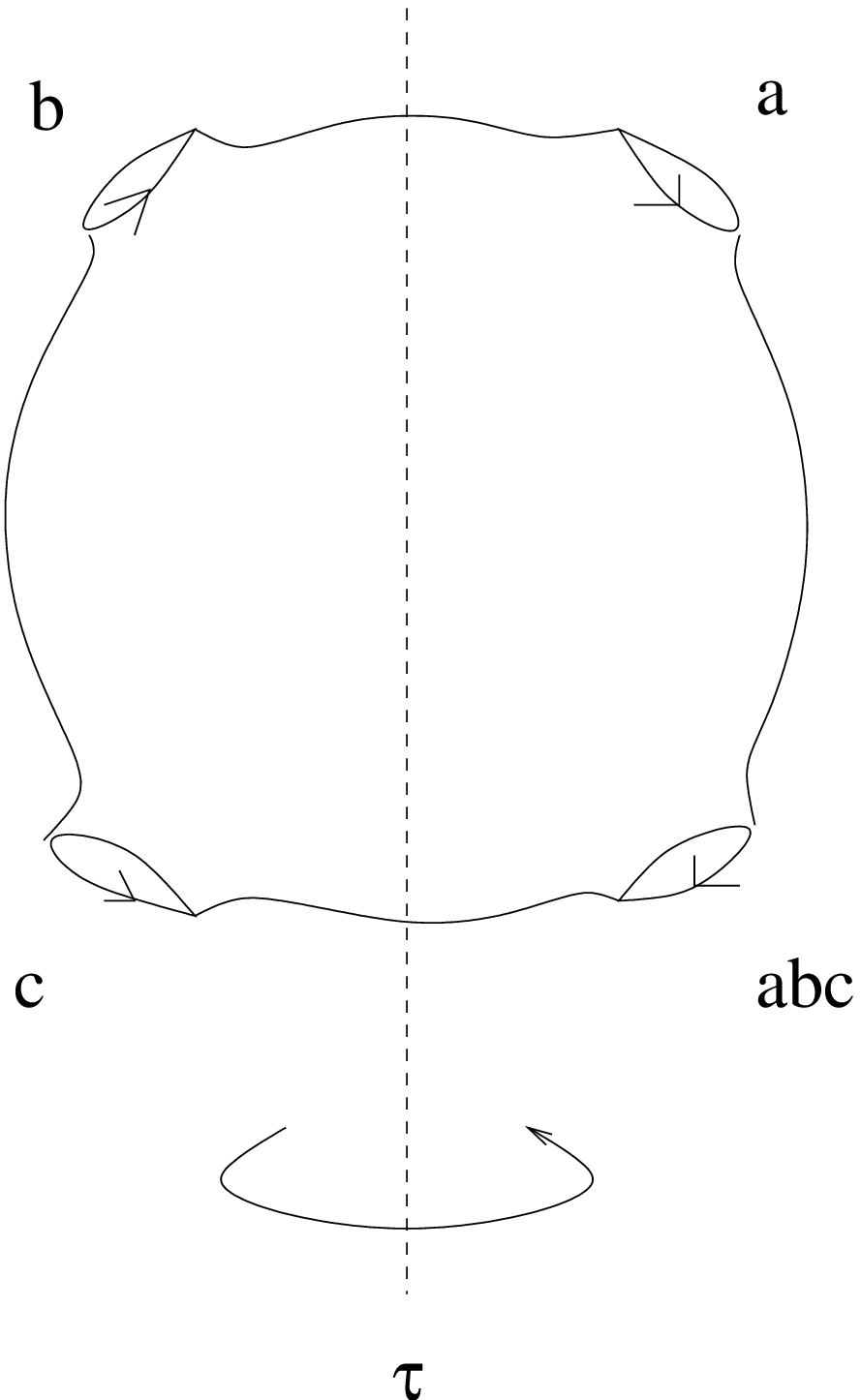}}
    \qquad\qquad
    \subfigure[$T_1$]{
      \includegraphics[height=3.2cm]{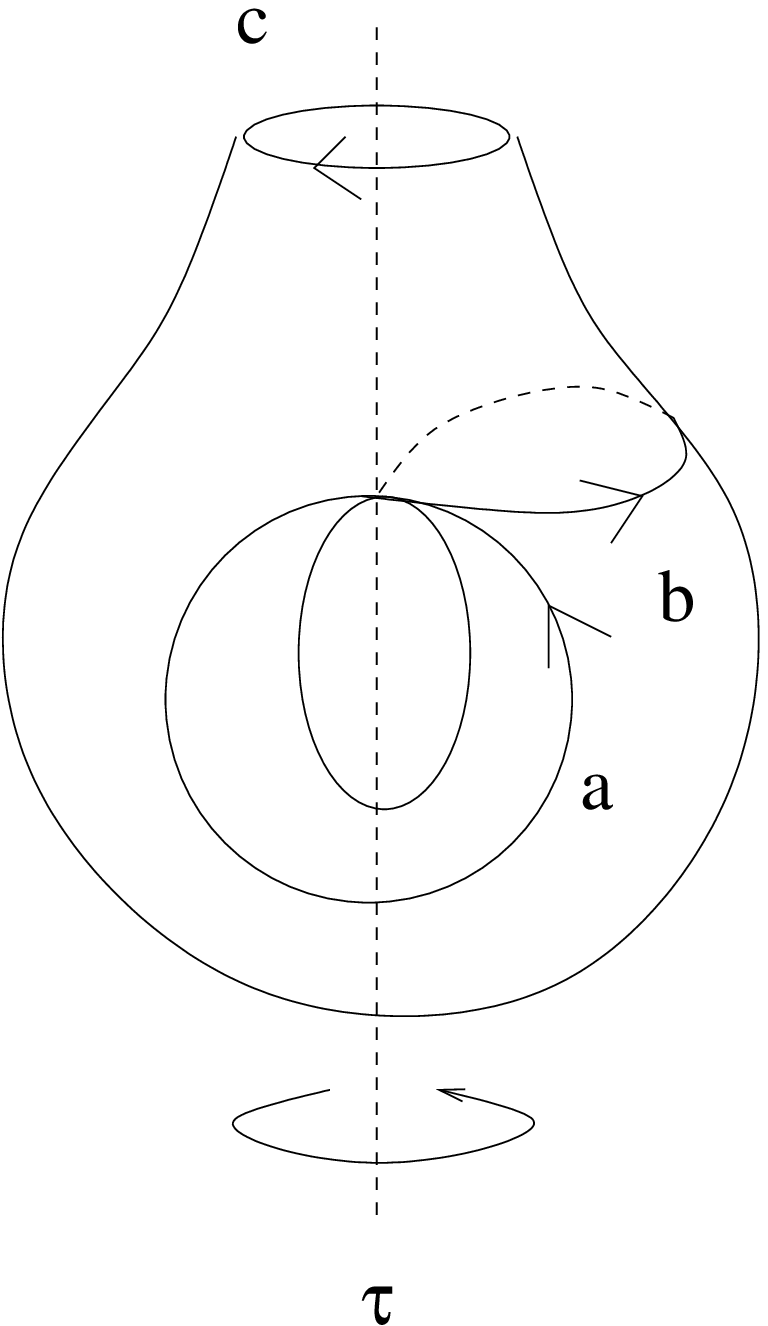}}
     \qquad
    \subfigure[$T_2$]{
      \includegraphics[height=3.2cm]{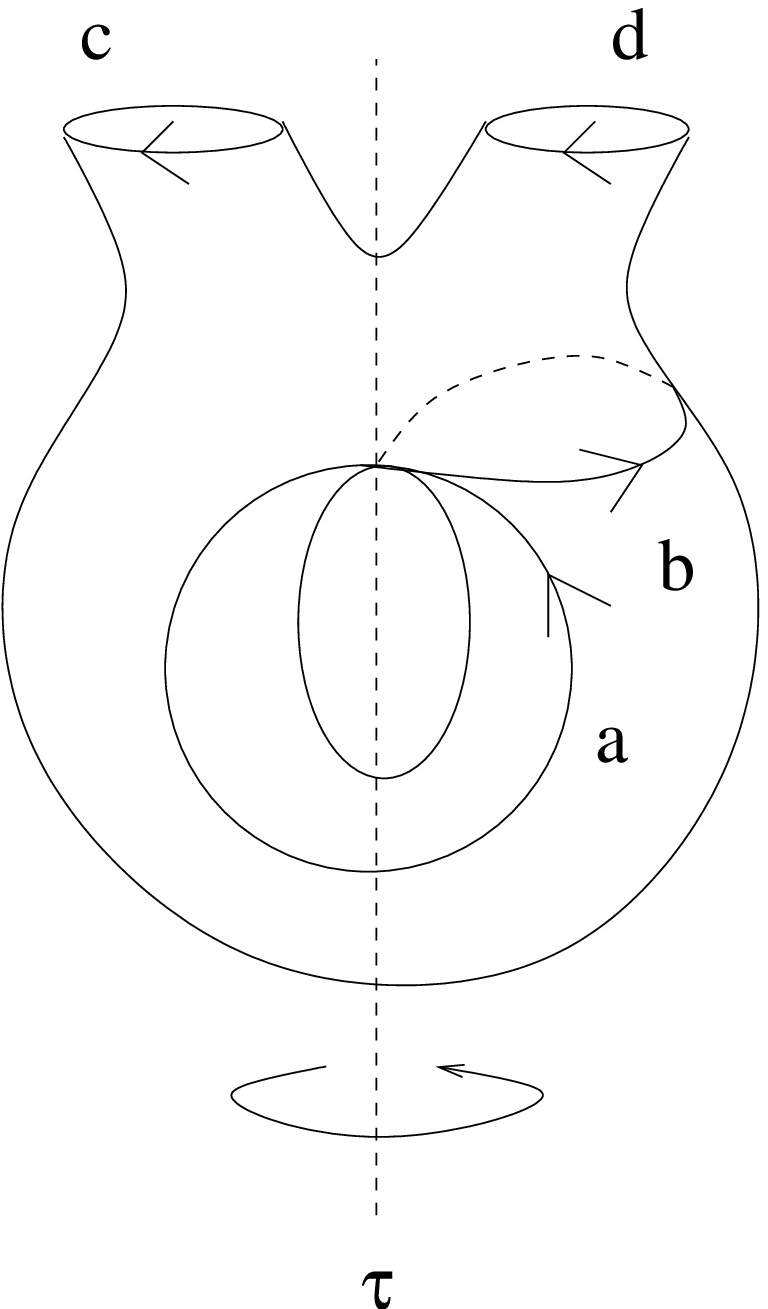}}
      \\
       \subfigure[$G_2$]{
      \includegraphics[height=2.5cm]{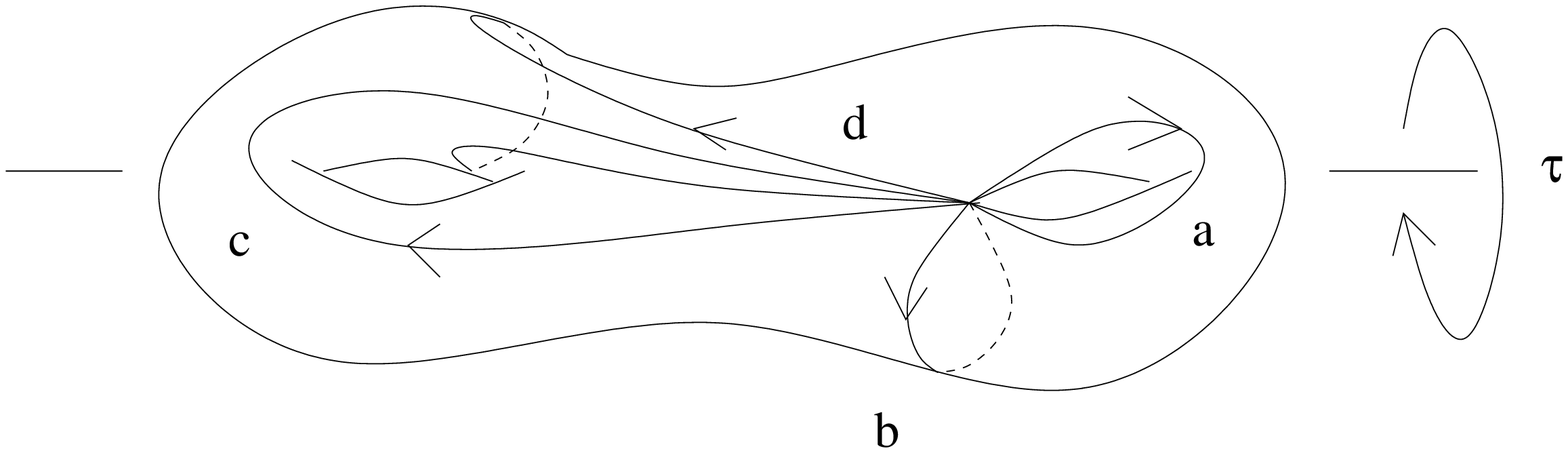}}
  \end{center}
  \caption{Ruberman's symmetric surfaces: a
      thrice-punctured sphere $S_3$, a four-punctured sphere $S_4$, a
      once-punctured torus $T_1$, a twice-punctured torus $T_2$ and a genus-two
      surface $G_2$ along with the shown involutions.} 
  \label{fig:mutsurf}
\end{figure}

We construct a birational map between certain subvarieties of the character
varieties of $M$ and $M^\tau$, which shows that in many cases the character
varieties are birationally equivalent. A subvariety $\T (M)$ in the
$\SL$--character variety $\X (M)$ will be defined, and the birational
equivalence is defined for subvarieties of $\T (M)$ which contain a dense set
of \emph{mutable} characters. All these notions descend to the
$\PSL$--character variety, and the objects are denoted by placing a
bar over the previous notation. 

\begin{pro} \label{mut: main pro}
Let $(S, \tau)$ be a mutation surface in an irreducible
3--manifold $M$, and let $C$ be an irreducible component of $\PT (M)$ which
contains the character of a representation whose restriction to 
$\im(\pi_1(S)\to\pi_1(M))$ is irreducible and has trivial centraliser.

Then $C$ is birationally equivalent to an irreducible component of
$\PT (M^\tau)$.
\end{pro}

\begin{pro} \label{mut: sep pro}
Let $(S, \tau)$ be a \emph{separating} mutation surface in an irreducible
3--manifold $M$, and let $C$ be an irreducible component of $\T (M)$ which
contains the character of a representation whose restriction to
$\im(\pi_1(S)\to \pi_1(M))$ is irreducible. 

Then $C$ is birationally equivalent to an irreducible component of
$\T (M^\tau)$.
\end{pro}

If $M$ admits a complete hyperbolic structure of finite volume, then there 
is a discrete and faithful representation of $\pi_1(M)$ into $\PSL$,
which lies on the \emph{Dehn surgery component} $\PX_0(M)$ 
of $\PX(M)$, and lifts to a component
$\X_0(M)$ of $\X(M)$. We now focus on the 
\emph{symmetric surfaces} shown in Figure \ref{fig:mutsurf}.

\begin{cor}[Hyperbolic knots] \label{mut: knots cor}
  Let $\k$ be a hyperbolic knot and $\k^\tau$ be a Conway mutant of $\k$.
  Then $\PX_0 (\k )$ and $\PX_0 (\k^\tau )$, as well as
  $\X_0 (\k )$ and $\X_0 (\k^\tau )$, are birationally equivalent.
  Moreover, the associated factors of the $A$--polynomials are equal.
\end{cor}

It is noted in \cite{tillus} that any Conway mutation of a knot can be
realised by at most two mutations along genus two surfaces.
The previous corollary is therefore a special case of the following:

\begin{cor}[Separating in hyperbolic] \label{mut: sep cor}
Let $(S, \tau)$ be a separating symmetric surface in a finite volume
hyperbolic 3--manifold $M$.
If $S$ is a twice--punctured torus or a genus two surface, 
then $\PX_0(M)$ and $\PX_0(M^\tau )$ as well as
$\X_0(M)$ and $\X_0(M^\tau )$ are birationally equivalent.
\end{cor}

The restriction to the two surfaces in the above corollary is necessary
in general. There are no separating incompressible and
$\partial$--incompressible thrice punctured spheres and once--punctured tori in
hyperbolic 3--manifolds, and it is easy to find examples of
Conway mutation on links with the property that only a proper subvariety in
each of the respective Dehn surgery components is contained in $\T (M)$ and
$\PT (M)$. If the surface is non--separating,
one can similarly find examples such that
mutation along a twice--punctured torus or a thrice punctured sphere does
not allow a general statement, which limits us to the following:

\begin{cor}[Non--separating in hyperbolic] \label{mut: non-sep cor}
Let $(S, \tau)$ be a non--separating symmetric surface in a finite
volume hyperbolic 3--manifold $M$.
If $S$ is a once--punctured torus or a genus two surface, then $\PX_0(M)$
and $\PX_0(M^\tau )$ are birationally equivalent.
\end{cor}

This corollary does not extend to $\SL$--Dehn surgery components in
general: mutation of the figure eight knot complement along the
fibre results in the associated sister manifold, and
the smooth projective models of their $\SL$--Dehn surgery components
are a torus and a sphere respectively.

The proofs of the above results are contained in
Section \ref{section:Tentatively mutable representations}.
Some of the ideas in the proofs are useful
in other settings; e.g.\thinspace
they produce examples of ``holes in the eigenvalue variety'' in \cite{tillus}.
The \emph{extension lemma} (see Lemma \ref{mut: extension lem}) can be used to
study the character variety of a 3--manifold by successively cutting along
non--separating surfaces.

In certain cases, analysis of the points where the
birational equivalence is not well--defined can be used to decide whether a
mutation surface is detected by an ideal point of the character variety. 
It is still an open problem whether every essential surface is detected by an
ideal point of the character variety. 
Necessary and sufficient conditions which have to be satisfied by a connected
surface are given in Section \ref{section:surfaces and ideal points},
and the birational equivalence is used to show that symmetric surfaces are
detected in the complements the Kinoshita--Terasaka knot and the figure eight
knot, as well as the so--called sister manifold of the latter.


\rk{Acknowledgements}
The contents of this paper forms part of the author's dissertation, and he
thanks his advisors Craig Hodgson and Walter Neumann for their
constant support and inspiration. He also thanks the referee and
Joan Birman for helpful comments on the previous version of this paper.
This work was supported by an
International Postgraduate Research Scholarship by the Commonwealth of
Australia Department of Education, Science and Training.


\section{Tentatively mutable representations}
\label{section:Tentatively mutable representations}

Our standard references for character varieties are \cite{cs, bozh}.
We recall some definitions and facts. The \emph{$\SL$--representation variety}
of a finitely generated group $\G$ is $\R (\G) = \Hom (\G, \SL )$.
Each $\rho \in \R (\G)$ defines a \emph{character} $\chi_\rho : \Gamma \to \C$
by $\chi_\rho (\gamma ) = \tr\rho (\gamma )$, and the set of characters $\X
(\G)$ is the \emph{$\SL$--character variety}. 
Both varieties are regarded as 
affine algebraic sets, and there is a regular map 
$\te : \R (\G) \to \X (\G)$. 
If $\G$ is the fundamental group of a topological space $M$, we write
$\R (M)$ and $\X (M)$ instead of $\R (\G)$ and $\X (\G)$ respectively.

A representation is \emph{irreducible} if
the only  subspaces of $\C^2$ invariant under its image are trivial.
Otherwise it is \emph{reducible}.
Irreducible representations are determined by characters up to 
inner automorphisms of $\SL$.
Let $\R^i (\G )$ denote the closure of the set of irreducible
representations, then the images
$\X^r(\G ) = \te (\Red (\G ))$ and
$\X^i(\G ) = \te (\R^i (\G ))$ are closed algebraic sets, and we have
$\X^r(\G ) \cup \X^i(\G ) = \X (\G )$.
The variety $\X^r(\G )$ is completely determined by the abelianisation of
$\G$.

There is a character variety arising from representations
into $\PSL$, and the relevant objects are
denoted by placing a bar over the previous notation. 
As with the $\SL$--character variety, the surjective
map $\pte : \PR (\G) \to \PX (\G)$ is
constant on conjugacy classes and if $\prho$ is an
irreducible 
representation, then $\pte^{-1}(\pte (\prho ))$ is the orbit of $\prho$
under conjugation. 
The natural map 
$\qe: \X (\G) \to \PX (\G)$ is
finite--to--one, but in general not onto. It is the quotient map
corresponding to the action of $\Hom (\Gamma, \Z_2)$ on $\X (M)$.
This action is not free in general.


\subsection{Tentatively mutable in $\mathbf{\SL}$}

Given a mutation surface $(S,\tau )$, we define a subvariety in $\R (S)$ by
\begin{equation*}
\R_\tau (S) = \{ \rho \in \R (S) \mid \tr\rho (\gamma ) = \tr\rho (\tau_*
\gamma) \thinspace \text{for all} \thinspace\thinspace \gamma \in \pi_1(S) \}.
\end{equation*}
This subvariety descends to the character variety, 
and we let $\X_\tau (S) = \te (\R_\tau (S))$.
In fact, $\tau$ induces a polynomial automorphism of $\X (S)$, and
$\X_\tau (S)$ is the set of its fixed points.
If $\R_{\tau}(S)$ contains an irreducible representation, then the
subvariety of reducible representations has positive codimension.
For the symmetric surfaces, one obtains the following result:

\begin{lem}{\rm\cite{tillus}}\qua\label{firstlemma}\label{codim}%
Let $(S,\tau)$ be a symmetric surface as described in Figure
  \ref{fig:mutsurf}. If $S=T_1$ or $S=G_2$, then $\R (S) = \R_\tau (S)$.
  Otherwise the character of $\rho \in\R(S)$ is invariant under $\tau$
  if and only if it satisfies the following equations:
  \begin{itemize}
    \item if $S = S_3$, $\pi_1(S_3)=\langle a,b\rangle$ then $\tr\rho (a) = \tr\rho (b)$,
    \item if $S = S_4$, $\pi_1(S_4)=\langle a,b,c\rangle$ then $\tr\rho (a) = \tr\rho (b)$
          and $\tr\rho (c) = \tr\rho (abc)$,
    \item if $S = T_2$, $\pi_1(T_2) = \langle a,b,c\rangle$ then $\tr\rho (c) = \tr\rho
          (c^{-1}[a,b])$.
  \end{itemize}
  
  The subvariety of reducible
  representations in $\R_\tau (S)$ has codimension one.
  Moreover, this property is preserved under $\te$.
\end{lem}

If $(S, \tau )$ is a mutation surface in a 3--manifold $M$, we call a
representation  $\rho\in\R (M)$ 
\emph{tentatively mutable with respect to $(S, \tau )$} if its character
restricted to $S$ is invariant under $\tau$. The set $\S (M)$ of
these representations is a subvariety of $\R (M)$. Let $\te (\S (M)) = \T
(M)$. The $S$--reducible characters form a closed set in $\T (M)$, which we
denote by $\F (M)$. Let the closure of $\T (M) - \F (M)$ in $\T (M)$ be $\M
(M)$.\footnote{The varieties $\mathfrak{S}$, $\mathfrak{T}$ and $\F$ of
  \cite{tillus} are here denoted by $\R_\tau$, $\T$ and $\F^i$ respectively.}
Then $\M (M) \subseteq \X^i(M)$ is the union of irreducible components of $\T
(M)$ which contain the character of a $S$--irreducible representation. 

In particular, if $(S,\tau)$ is one of the symmetric surfaces $(T_1, \tau)$ or
$(G_2,\tau)$, we have $\T (M) = \X (M)$, so $\X^r(M)\subseteq \F (M)$, and the
same is true for $M^\tau$. In general, it is not true that $\T (M) = \X (M)$
implies $\T (M^\tau ) = \X (M^\tau )$.


\subsection{Tentatively mutable in $\mathbf{\PSL}$}

Given a mutation surface $(S,\tau )$, we define a subvariety in $\PR (S)$ by
\begin{equation*}
\PR_\tau (S) = \{ \prho \in \PR (S) \mid \chi_{\prho} = \chi_{\prho\tau_*} \}.
\end{equation*}
This subvariety descends to the character variety, 
and we let $\PX_\tau (S) = \pte (\PR_\tau (S))$.
For the symmetric surfaces, we have the following lemma.%
\footnote{The author thanks Steven Boyer for
pointing out that an earlier version of this lemma was incorrect.}

\begin{lem}\label{psl-firstlemma}
  Let $(S,\tau)$ be a symmetric surface as described in Figure
  \ref{fig:mutsurf} and $\prho$ be a $\PSL$--representation
  of $\pi_1(S)$. If $S$ is one of the surfaces with boundary, then
  $\chi_{\prho} = \chi_{\prho\tau_*}$ if and only if there is a lift
  $\rho$ of $\prho$ such that $\chi_{\rho} = \chi_{\rho\tau_*}$. 
\newpage
  If $S=G_2$, then $\prho \in \PR_\tau(G_2)$ either if $\prho$ lifts to a
  $\SL$--representation, or if
{\small
\begin{equation*}
    (\tr\prho(ad^{-1}))^2 =
   ( \tr\prho(bc^{-1}))^2 =
       ( \tr\prho(abd^{-1}))^2 =
            (\tr\prho(b^{-1}cd))^2 =
                (\tr\prho(acd))^2 = 0.
\end{equation*}
}\rm
\end{lem}

\begin{proof}
We use Lemma 3.1 of \cite{bozh} throughout this proof, and assume
familiarity with the notation used there; in particular, $\epsilon$ denotes
a homomorphism into the group $\{ \pm 1\}$.

Let $\prho \in \PR (S)$, and assume that there is a lift $\rho$ of $\prho$
such that $\rho \in \R_\tau (S)$. It then follows
that $\prho \in \PR_\tau (S)$, by choosing $\rho' = \rho\tau_*$
and $\epsilon = \identity$.

Since the fundamental groups of the surfaces with boundary are free, every
$\PSL$--representation of these surfaces lifts to a $\SL$-representation. We
now verify the statement of the lemma for these surfaces.

{\bf Case $S = T_1$}\quad Since $\PR (T_1) = \qe (\R(T_1))$ and
$\R(T_1)=\R_\tau (T_1)$, there is nothing to prove. In particular, we have
$\PR(T_1)=\PR_\tau (T_1)$.

{\bf Case $S=T_2$}\quad Let $\prho \in \PR_\tau (S)$, and $\rho \in \R(S)$ be a
lift of $\prho$. We have
$\chi_{\prho} = \chi_{\prho\tau_*}$ if and only if there is
$\epsilon \in \Hom(\pi_1(S), \{\pm 1\})$ such that
$\epsilon \chi_{\rho} = \chi_{\rho\tau_*}$.
Now $\epsilon(a) \tr\rho (a) = \tr\rho\tau(a) = \tr\rho(a^{-1})$ forces
$\epsilon (a) =1$. Similarly,
$\epsilon(b) \tr\rho (b) = \tr\rho\tau(b) = \tr\rho(ab^{-1}a^{-1})$ forces
$\epsilon (b) =1$.
Then $\epsilon(bc) \tr\rho (bc) = \tr\rho\tau(bc) = \tr\rho(bc)$ yields
$\epsilon(c)=\epsilon(b)$. Thus, $\epsilon = \identity$, and the claim
follows.

{\bf Case $S=S_3$}\quad Let $\prho \in \PR_\tau (S)$, and $\rho \in \R(S)$ be a
lift of $\prho$. Since $\tau(a)=b$, we have
$(\tr\prho(a))^2 = (\tr\prho(b))^2$. If $\tr\rho(a) = \tr\rho(b)$, then
$\rho \in \R_\tau(S)$. If $\tr\rho(a) = -\tr\rho(b)$, then define
$\sigma(a)=\rho(a)$ and $\sigma(b)=-\rho(b)$. Then $\sigma$ is a lift of
$\prho$ and $\sigma \in \R_\tau(S)$. This completes the proof in this
case.

{\bf Case $S=S_4$}\quad Let $\prho \in \PR_\tau (S)$, and $\rho \in \R(S)$ be a
lift of $\prho$. We have
$\chi_{\prho} = \chi_{\prho\tau_*}$ if and only if there is a homomorphism
$\epsilon \in \Hom(\pi_1(S), \{\pm 1\})$ such that
$\epsilon \chi_{\rho} = \chi_{\rho\tau_*}$.
As above, considering the action of $\tau$ yields
$\epsilon(a)=\epsilon(b)=\epsilon(c)$. If $\epsilon$ is trivial, then
$\rho \in \R_\tau(S)$. Otherwise, the character of the lift $\sigma$ defined
by $\sigma(a) = \rho(a)$, $\sigma(b)=\rho(b)$ and $\sigma(c) = - \rho(c)$ is
invariant under $\tau$.

Now consider $S=G_2$. It follows from Theorem 5.1 in \cite{g}
that $\PX(G_2)$ has two topological components with the property that every
representation in one of the components lifts to $\SL$, and every
representation in the other does not. We only
have to consider the latter component since $\R (G_2) = \R_\tau(G_2)$.

Assume that $\prho$ is a $\PSL$--representation of $G_2$ with representative
matrices $A,B,C,D$ for the generators $a,b,c,d$,
such that $[A,B][C,D]=-E$. Then $\prho$ does not
lift to $\SL$. Now assume that $\prho \in \PR_\tau (S)$, and define a
representation $\rho \in \R(\F_4)$ by $\rho(\alpha)=A$,
$\rho(\beta)=B$, $\rho(\gamma)=C$ and $\rho(\delta)=D$.
By assumption, there is $\epsilon \in \Hom(\F_4, \{\pm 1\})$ such that
$\epsilon \chi_\rho = \chi_{\rho\tau}$, where $\rho\tau$ is defined by
\begin{align*}
    \rho\tau (\alpha) &= A^{-1}, &&& \rho\tau (\gamma) &= (B^{-1}CD) C^{-1}
    (B^{-1}CD)^{-1},                   \\
    \rho\tau (\beta)  &= AB^{-1}A^{-1}, &&&  \rho\tau (\delta) &= (B^{-1}C)
    D^{-1} (B^{-1}C)^{-1}.          
\end{align*}
Then
$\epsilon(a) \tr A = \epsilon(a) \tr\rho (\alpha)
= \tr\rho\tau (\alpha) = \tr A^{-1}$
forces $\epsilon(a)=1$. We similarly obtain
$1 = \epsilon(b)= \epsilon(c) = \epsilon(d)$. But then
$\epsilon =\identity$, and we have
{\small{\begin{align*}
       \tr (AD)
      =& \epsilon(\alpha\delta) \tr\rho(\alpha\delta)
      = \tr\rho\tau (\alpha\delta)\\
      =& \tr (A^{-1}(B^{-1}CD^{-1}C^{-1}B))
            &&\mid\text{by definition of $\rho\tau$}\\
      =& \tr(A)\tr(D) - \tr (AB^{-1}CD^{-1}C^{-1}B)
            &&\mid\text{by $\tr X^{-1}Y = \tr X \tr Y - \tr XY$}\\
      =& \tr (A)\tr\rho (D) + \tr (AB^{-1}D^{-1}ABA^{-1})
            &&\mid\text{by $[A,B][C,D]=-E$}\\
      =& \tr(A)\tr (D) + \tr (D^{-1}A) 
            &&\\
      =& \tr(AD) + 2 \tr (D^{-1}A) 
            &&\mid\text{by $\tr X^{-1}Y = \tr X \tr Y - \tr XY$}. 
\end{align*}
}\rm}\rm
Thus, $\tr(AD^{-1})=0$, and therefore $(\tr\prho(ad^{-1}))^2=0$. The other
stated trace identities follow similarly. 
This completes the proof of the lemma.
\end{proof}

If $\prho \in \PR_\tau (S)$ is an irreducible representation, then there
exists an element $\overline{X} \in \PSL$, such that 
$\prho = \overline{X}^{-1} \prho\tau \overline{X}$. 
The centraliser of an element $\overline{Y}$ in $\PSL$
is the quotient of its centraliser in $\SL$ unless
$(\tr \overline{Y})^{2} = 0$. 
Thus, if $\G$ is a finitely generated group, then
the centraliser of an irreducible representation $\prho \in\PR(\G)$ is trivial
if $(\tr\prho (\gamma_i))^2 \ne 0$ for all 
generators $\gamma_i$ of $\G$. 
Let $\F_2$ be the free group in two elements $\langle\alpha, \beta \rangle$.
We have $\X(\F_2) \cong \C^3$, and the map $\PX (\F_2) \to \C^3$ given by
$\chi_{\prho} \to ( (\tr\prho(\alpha))^2, (\tr\prho(\beta))^2,
(\tr\prho(\alpha\beta))^2)$ is a $2:1$ covering map.

\begin{lem} \label{mut: free}
Consider the above two--to--one parameterisation of the $\PSL$--character
variety of $\F_2 = \langle\alpha, \beta \rangle$ by the points
$((\tr\prho \alpha)^2, (\tr\prho \beta)^2,(\tr\prho \alpha\beta)^2)$ in
$\C^3$. Then the set of irreducible representations with non--trivial
centraliser is contained in the union of the three coordinate axes. 
\end{lem}

\begin{proof}
Assume that $\prho$ is an irreducible representation of $\F_2$ with
non--trivial centraliser in $\PSL$. According to the above discussion
at least one of $(\tr\prho \alpha)^2$ or $(\tr\prho \beta)^2$ is equal to
zero. Assume that $(\tr\prho \alpha)^2=0$. Direct calculation shows that the
centraliser of $\prho(\F_2)$ is non--trivial if and only if 
$(\tr\prho \beta)^2=0$ or 
$(\tr\prho \alpha\beta)^2=0$. If both are equal to zero, then the image of
$\prho$ is a Kleinian four group in $\PSL$ and equal to its centraliser, and
if one of the traces is not equal to zero, then the centraliser has order
equal to two. 
\end{proof}

It follows that if $(S, \tau)$ is a mutation surface and $\PR_\tau (S)$
contains an irreducible representation with trivial centraliser, then the set
of reducible representations and the set of
irreducible representations with non--trivial centraliser are contained in
subvarieties of positive codimension. In particular:

\begin{lem} \label{mut: psl-codim}
  Let $(S,\tau )$ be a symmetric surface. The set of reducible representations
  in $\PR_\tau (S)$
  and the set of representations in $\PR_\tau (S)$
  with non--trivial centralisers are contained in a finite union of
  subvarieties, each of which has codimension one.
  Moreover, this property is preserved under $\pte$.
\end{lem}

\begin{proof}
  The subvariety of reducible representations has
  codimension one since the proof of Lemma \ref{codim} 
  (Lemma 2.1.3 in \cite{tillus}) applies again.
  The set of irreducible representations with non--trivial centralisers are
  contained in a union of subvarieties each of which is defined by stating
  that two coordinates are equal to zero. Each of these subvarieties is
  easily observed to have codimension at least one in $\PR_\tau(S)$ for each
  of the symmetric surfaces.
\end{proof}

We can now define $\PM (M)$ to be the union of the irreducible components of 
$\PT (M)$ which contain the character of an $S$--irreducible 
representation such that
the image of $\im(\pi_1(S)\to\pi_1(M))$ has
trivial centraliser.


\subsection{Extension lemma}

Let $A$ be a finitely generated group and $\varphi :A_1 \to A_2$ be an
isomorphism between finitely generated subgroups of $A$. Define
\begin{equation*}
  \PR_\varphi (A) := \{ \prho \in \PR (A) \mid \chi_{\prho}|_{A_1} = 
  \chi_{\prho\varphi}|_{A_1} \},
\end{equation*}
and $\pte (\PR_\varphi(A) )= \PX_\varphi(A)$. Let
$\G=\langle A,k\mid k^{-1}ak=\varphi (a)\thinspace\forall a\in A_1\rangle$
be a HNN--extension of $A$. Assume that
$\prho \in \PR_\varphi (A)$ has the property that
$\prho|_{A_1}$ is irreducible with
trivial centraliser. Then there exists a unique $\prho' \in \PR(\G)$ 
such that $\prho'|_{A} = \prho$: the assignment 
$\prho'(k)$ is the unique element of $\PSL$ which conjugates
$\prho$ to $\prho\varphi$.

\begin{lem} \label{mut: extension lem}
Let
$\G$ and $\PX_\varphi (A)$ be as defined above. Let $V$ be an irreducible
component of $\PX (\G )$ containing the character of a
representation which restricted to $A_{1}$ is irreducible and has trivial
centraliser. Then the restriction map $r : \PX (\G) \to \PX_\varphi (A)$ is a
birational equivalence between $V$ and $\overline{r(V)}$. 
\end{lem}

\begin{proof}
The restriction map is a polynomial map, and hence $W:=\overline{r(V)}$ is
an irreducible component of $\PX_\varphi (A)$. If follows from Lemma 
\ref{mut: free}, Lemma 4.1 of \cite{bozh} and the fact that irreducible
representations with the same character are equivalent, that the above
construction of $\PSL$--representations of $\G$ from $\PSL$--representations
of $A$ is a well--defined
1--1 correspondence of $\PSL$--characters in $V$ and $W$
apart from a subvariety of codimension at least one. Thus, $r$ has degree
one and is therefore a birational isomorphism onto its image.
\end{proof}


\subsection{Proofs of the main results}


\rk{Proof of Proposition \ref{mut: sep pro}}
The following construction is taken from \cite{a-poly1}.
Given a separating mutation surface $(S, \tau)$, we obtain a decomposition 
\begin{equation*}
  \pi_1(M) \cong \pi_1(M_{-}) \star_{\pi_1(S)} \pi_1(M_{+}).
\end{equation*}
The varieties $\R (M)$ and $\R (M^\tau)$ can be viewed as a subsets of 
$\R (M_- )\times \R (M_+ )$, and the inclusion map is the
restriction to the respective subgroups.
Let $\rho \in \S (M)$ be an $S$--irreducible  representation.
Since $\rho_-\tau$ is equivalent to $\rho_-$ on $\pi_1(S)$, there is an
element $X \in \SL$ 
such that $\rho_- = X^{-1} \rho_-\tau X$ on $\pi_1(S)$, and
$X$ is defined up to sign.
We can now define a representation $\rho^\tau$ of $M^\tau$
as follows: Let
$\rho_+^\tau = \rho_+$ on $\pi_1(M_{+})$ and
$\rho_-^\tau = X^{-1} \rho_-\tau X$
on $\pi_1(M_{-})$. Then $\rho^\tau = (\rho_-^\tau ,\rho_+) \in \R (M^\tau)$
is well defined, since both definitions agree on the amalgamating subgroup,
and the map $\rho \to \rho^\tau$ only depends upon the inner automorphism
induced by $X$.
Both $\rho$ and $\rho^\tau$
are irreducible and $\rho^\tau \in \S (M^\tau )$.
It is shown in \cite{tillus} that this construction yields an isomorphism
$\map : \M (M) \to \M (M^\tau )$ defined everywhere apart from the subvariety
$\F^i (M)$ of characters of irreducible representations which are reducible
on $\pi_1(S)$. Moreover, it is shown on pages 567-568 of \cite{tillus} that
$\map$ is a birational equivalence between irreducible components 
(since they contain a $S$--irreducible character). \qed


\rk{Proof of Proposition \ref{mut: main pro}}
Assume that $S$ is separating.
The previous construction of representations also works for
projective representations with trivial centraliser, and the argument in the
above mentioned proof goes through if one uses Lemma 4.1 of \cite{bozh}
instead of Proposition 1.1.1 of \cite{cs}.

Thus, let $S$ be a non--separating mutation surface.
The boundary of $M-S$
contains two copies $S_+$ and $S_-$ of $S$.
Let $A = \im(\pi_1(M-S)\to \pi_1(M))$, $A_1=\im(\pi_1(S_+)\to \pi_1(M))\le
A$ and $A_2 = \im(\pi_1(S_-)\to\pi_1(M)) \le A$. Then 
$\pi_1(M)$ is an HNN--extension of
$A$ by some $k \in \pi_1(M)$ across $A_1$ and $A_2$:
\begin{equation*}
  \pi_1(M) = \langle A,k \mid k^{-1}A_1k = A_2\rangle.
\end{equation*}
The action of $k$ is determined by the gluing map $S_+ \to S_-$, and the
mutation changes the gluing map by $\tau$.
We thus obtain a presentation of $\pi_{1}(M^\tau)$:
\begin{equation*}
  \pi_1(M^\tau) = \langle A,k \mid k^{-1}\tau(A_1)k = A_2\rangle.
\end{equation*}
Let $\prho$ be a $\PSL$--representation of $M$.
Note that $\prho (k)$ is only determined up to the centraliser of $\prho
(A_1)$. 
Assume that $\prho$ is tentatively mutable and $\prho (A_1)$ is irreducible
and has trivial centraliser. 
Then $\prho (k)$ is uniquely determined by $\prho (A_1)$ and the gluing map.
Furthermore, $\prho (a)$ is conjugate to $\prho \tau (a)$ via
some uniquely determined $X \in \PSL$ for all $a\in A_1$. It follows that
$\prho (k)^{-1}\prho (a_1) \prho (k) = \prho (a_2)$ is equivalent to
$\prho (k)^{-1} X^{-1}\prho \tau(a_1) X \prho (k) = \prho (a_2)$.
Define a representation $\prho^\tau \in \PR(M^\tau)$ by
$\prho^\tau (a) := \prho (a)$ for all $a \in A$, and
$\prho^\tau (k) = X \prho (k)$.
Denote the corresponding map by $\pmap$.
Since we can define an inverse map, we have a
natural 1--1 correspondence of $A_1$--irreducible
representations (with non--trivial centraliser on $A_1$)
in $\PS (M)$ and $\PS (M^\tau)$.
Moreover, this map is well--defined for equivalence classes of
representations, and hence for the corresponding characters in 
$\PM (M)$ and  $\PM (M^\tau)$.

Let $C$ be an irreducible component of $\PM (M)$, i.e.\thinspace a component 
of $\PT (M)$ which
contains the character of a $S$--irreducible $\PSL$--representation such that
the image of $\pi_1(S)$ has trivial centraliser.
By definition, the restriction maps
$r : \PM (M) \to \PX_\varphi (A)$ and
$r^\tau : \PM (M^\tau) \to \PX_{\varphi\tau} (A)$ have range in a subvariety
of $\PX_\varphi (A) \cap \PX_{\varphi\tau} (A)$.
The construction of $\pmap$ gives $r ( \chi ) = r^\tau (\pmap \chi)$,
whenever applicable. Since $\pmap$ is defined on a
dense subset of $C$, Lemma \ref{mut: extension lem}
implies that it is the composition
$(r^\tau)^{-1}\circ r  $. \qed


\rk{Proof of Corollary \ref{mut: sep cor}}
Assume that $M$ is a finite volume hyperbolic 3--manifold and
$S$ is a separating symmetric surface and either $T_2$ or $G_2$. 
If $\X_0(M) \subseteq \T(M)$, then 
$\PX_0(M) \subseteq \PT(M)$, since $\qe (\X_0) = \PX_0$ and  
Lemma \ref{psl-firstlemma} applies.
The two boundary components of any separating incompressible $T_2$ have to
lie on the same boundary component of $M$, hence Lemma \ref{firstlemma}
implies that $\T (M) = \X (M)$.
Since $\rho_0$ is torsion free and $S$--irreducible, both
$\PX_0(M)$ and $\X_0(M)$ satisfy the hypotheses of Propositions
\ref{mut: main pro} and \ref{mut: sep pro} (where applicable).
It follows from \cite{r} that the birational equivalence takes the complete
representation of $M$ to the complete representation of $M^\tau$, 
and hence it restricts to
a birational equivalence between the two Dehn surgery components.\qed


\rk{Proof of Corollary \ref{mut: non-sep cor}}
Since $S=T_1$ or $S=G_2$, we have $\X_0(M) \subseteq \T(M)$ and
$\X_0(M^\tau) \subseteq \T(M^\tau)$. The same arguments as in
the proof of Corollary \ref{mut: sep cor} now yield the conclusion. \qed


\rk{Remark}
The proofs of Propositions \ref{mut: main pro} and \ref{mut: sep pro} show
that we have birational equivalences $\map : \M (M) \to \M (M^\tau )$ and
$\pmap : \PM (M) \to \PM (M^\tau )$. Since every knot group abelianises to
$\Z$, this in particular implies:

\begin{pro} 
  Let $\k$ and $\k^\tau$ be Conway mutant knots.
  If every component of $\X (\k)$ and $\X (\k^\tau )$
  which contains the character of an irreducible representation contains the
  character of a $S$--irreducible representation, then
  $\X (\k)$ and $\X (\k^\tau )$ are birationally equivalent.
\end{pro}


\section{Surfaces and ideal points}
\label{section:surfaces and ideal points}

We build on the construction by Culler and Shalen \cite{sh1, bozh} to give a
method to determine whether a connected essential surface is associated to an
ideal point. This method is then applied to two pairs of mutative manifolds in
conjunction with the respective birational equivalences.


\subsection{Surface associated to the action}
\label{Surface associated to the action}

Let $M$ be an orientable, irreducible 3--manifold, and assume that $\tree_v$
is Serre's tree 
associated to an ideal point $\xi$ of a curve $C$ in $\X (M)$ or $\PX (M)$.
A surface \emph{associated to the action} of $\pi_1(M)$ on $\tree_v$
is defined by Culler and Shalen using a construction due to Stallings.
If the given manifold is not compact, replace it by a compact core.
Choose a triangulation of $M$ and give the universal cover $\tilde{M}$ the
induced triangulation. One can then construct a
simplicial, $\pi_1(M)$--equivariant map $f$ from $\tilde{M}$ to $\tree_v$. The
inverse image of midpoints of edges is a surface in $\tilde{M}$ which
descends to a non--empty, 2--sided surface $S$ in $M$.
The map $f$ is changed by a homotopy (if necessary) so that $S$ is
incompressible and has no boundary parallel or sphere components.
We then say that $S$ is \emph{essential}.
The associated surface $S$ depends upon the choice of triangulation of $M$
and the choice of the map $f$.
An associated surface often contains finitely many parallel copies of one of
its components. They are somewhat redundant,
and we implicitly discard them, whilst we still call the resulting surface
associated.


\subsection{Surface detected by an ideal point}
\label{Surface detected by ideal point}

We now describe associated surfaces satisfying certain non--triviality
conditions. An essential surface $S$ in $M$ gives rise to a graph of groups
decomposition of $\pi_1(M)$. Let $t_{1},...,t_{k}$ be the generators of the
fundamental group of the graph of groups arising from $HNN$--extensions.
Let $M_{1},...,M_{m}$ be the components of
$M-S$, let $\tree_S$ be the dual tree
to $\tilde{S}$ in $\tilde{M}$ and $\graph_S$ be the dual graph to $S$ in $M$.
For each component $M_i$ of $M-S$, fix a representative $\G_i$ of the
conjugacy class of $\im (\pi_1(M_i) \to \pi_1(M))$
as follows. Let $\tree' \subset \tree_S$ be a tree of representatives, 
i.e.\thinspace a
lift of a maximal tree in $\mathcal{G}_S$, and let
$\{ s_1,{\ldots} ,s_m\}$ be the vertices of $\tree'$, labelled such that $s_i$
maps to $M_i$ under the composite mapping $\tree_S \to \graph_S \to M$. Then
let $\G_i$ be the stabiliser of $s_i$.

Assume that $S$ does not contain parallel copies of one of its components.
Then $S$ is \emph{detected by an ideal point} of the character
variety with Serre tree $\tree_v$ if 
\renewcommand{\labelenumi}{\bf S\theenumi}
\begin{enumerate}
\item \label{S1} every vertex stabiliser of the action on $\tree_S$ is included in
a vertex stabiliser of the action on $\tree_v$,
\item \label{S2} every edge stabiliser of the action on $\tree_S$ is included in
an edge stabiliser of the action on $\tree_v$,
\item \label{S3} if $M_i$ and $M_j$, where $i \ne j$, are identified along a component
of $S$, then there are elements $\gamma_i \in \G_i$ and
$\gamma_j \in \G_j$ such that $\gamma_i\gamma_j$ acts as a loxodromic
on $\tree_v$,
\item \label{S4} each of the generators $t_i$ can be chosen to act as a
loxodromic on $\tree_v$. 
\end{enumerate}

\begin{lem} \label{lem:detected}
Let $M$ be an orientable, irreducible 3--manifold.
An essential surface $S$ in $M$ which is detected by an ideal point $\xi$ of a
curve $C$ in $\X (M)$ is associated to the action of $\pi_1(M)$ on the Serre
tree $\tree_v$. 
\end{lem}

\begin{proof}
Choose a sufficiently fine triangulation of $M$ such that the $0$--skeleton  
of the triangulation is disjoint from $S$, and such that the intersection of
any edge in the triangulation with $S$ consists of at most one point. Give
$\tilde{M}$ the induced triangulation. We may assume that the retraction
$\tilde{M} \to \tree_S$ is simplicial, 
and we now define a map $\tree_S \to \tree_v$.

The vertices $\{s_1,{\ldots} ,s_m\}$ of the tree of representatives
are a complete set of orbit representatives
for the action of $\pi_1(M)$ on the 0--skeleton of $\tree_S$.
Condition S\ref{S3} implies that we may choose vertices 
$\{ v_1,{\ldots} ,v_m\}$ of $\tree_v$ such that $v_i$ is stabilised by
$\G_i$, and if $M_i \ne M_j$, then $v_i \ne v_j$.
Define a map $f^0$ between the 0--skeleta of $\tree_S$ and $\tree_v$ as
follows. Let $f^0 (s_i) = v_i$. For each other vertex $s$ of $\tree_S$ there
exists $\gamma \in \pi_1(M)$ such that $\gamma s_i = s$ for some $i$. Then
let $f^0(s) = \gamma f^0(s_i)$. We thus obtain a
$\pi_1(M)$--equivariant map from $\tree_S^0 \to \tree_v^0$, which
extends uniquely to a map $f^1 : \tree_S \to \tree_v$, since the image of each
edge is determined by the images of its endpoints. Since
$v_i \ne v_j$ for $i \ne j$, and since each $t_k$ acts as a loxodromic on
$\tree_v$, the image of each edge of $\tree_S$ is a path of length greater or
equal to one in $\tree_v$.

If $f^1$ is not simplicial, then there is a subdivison of $\tree_S$ giving
a tree $\tree_{S'}$ and a $\pi_1(M)$--equivariant, simplicial map
$f : \tree_{S'} \to \tree_v$. There is a surface $S'$ in $M$ which is obtained
from $S$ by adding parallel copies of components such that $\tree_{S'}$ is the
dual tree of $\tilde{S}'$.

As before,
choose a sufficiently fine triangulation of $M$ such that the $0$--skeleton
of the triangulation is disjoint from $S'$, and such that the intersection of
any edge in the triangulation with $S'$ consists of at most one point, and
give $\tilde{M}$ the induced triangulation. The composite map
$\tilde{M} \to \tree_{S'} \to \tree_v$ is $\pi_1(M)$--equivariant and
simplicial, and the inverse image of midpoints of edges descends to the
surface $S'$ in $M$. Thus, $S'$ is associated to the action of $\pi_1(M)$ on
$\tree_v$.
\end{proof}

We now wish to decide whether a given essential surface $S$ in
$M$ is detected by an ideal point of a curve in $\X (M)$.
Denote the components of $M-S$ by $M_1, {\ldots} ,M_m$.
If $S$ is detected by an ideal point, then the limiting character restricted
to each $M_i$ is finite.
There is a natural map from $\X (M)$ to
$\X (M_1) \times {\ldots} \times \X (M_m)$ by restricting to
the respective subgroups.
Splittings along
$S$ which are detected by ideal points of curves in $\X (M)$ correspond to
points $(\chi_1,{\ldots} ,\chi_m)$ in the cartesian product satisfying the
following necessary conditions:  
\renewcommand{\labelenumi}{\bf C\theenumi}
\begin{enumerate}
\item \label{C1} $\chi_i \in \X (M_i)$ is finite for each $i=1,{\ldots} ,m$.
\item \label{C2} For each component of $S$, let  $\varphi : S^+ \to S^-$ be
  the gluing 
  map between its two copies arising from the splitting,
  and assume that $S^+ \subset \partial M_i$ and $S^- \subset \partial M_j$,
  where $i$ and $j$ are not necessarily distinct.
  Denote the homomorphism induced by $\varphi$ on fundamental group by
  $\varphi_*$. Then for each $\gamma \in \im (\pi_1(S^+) \to \pi_1(M_i))$,
  $\chi_i (\gamma ) = \chi_j (\varphi_* \gamma )$.
\item \label{C3} For each $i=1,{\ldots} ,m$, the restriction of $\chi_i$
      to any component of $S$ in $\partial M_i$ is reducible.
\item \label{C4} There is an ideal point $\xi$ of a curve $C$ in $\X (M)$
  and a  
  connected open
  neighbourhood $U$ of $\xi$ on $C$ such that
  the image of $U$ under the map to the cartesian product contains an open
  neighbourhood of $(\chi_1,{\ldots} ,\chi_m)$ on a curve in 
  $\X (M_1) \times {\ldots} \times \X (M_m)$, but not 
  $(\chi_1,{\ldots} ,\chi_m)$ itself. 
\end{enumerate}
The first condition implies that $\im(\pi_1(M_i) \to \pi_1(M))$ is contained
in a vertex stabiliser for each $i=1,{\ldots} ,m$. The second defines a
subvariety of the cartesian product containing the image of $\X (M)$ under the
restriction map. Condition
C\ref{C3} must be satisfied since it is shown in \cite{sh1} that the limiting
representation of every component of an associated surface is reducible.
The last condition implies that the action of $\pi_1(M)$ on Serre's tree is
non--trivial. 

\begin{lem} \label{associating lem}
Let $S$ be a connected essential surface in an orientable, irreducible
3-manifold $M$. Then $S$
is associated to an ideal point of the character variety of $M$ if and only if
there are points in the cartesian product of the character varieties of the
components of $M-S$ satisfying conditions C\ref{C1}--C\ref{C4}.
\end{lem}

\begin{proof}
We need to show that the conditions are sufficient. Assume that $S$ is
non--separating. Let $A = \im(\pi_1(M-S) \to \pi_1(M))$, and denote the
subgroups of $A$ corresponding to the two copies of $S$ in $\partial (M-S)$
by $A_1$ and $A_2$. Then $\pi_1(M)$ is an HNN--extension of $A$ 
by some $t \in \pi_1(M)$ across
$A_1$ and $A_2$. We may assume that $t^{-1} A_1 t = A_2$.

Let $\xi$ be the ideal point provided by C\ref{C4}, and denote Serre's tree
associated to $\xi$ by $\tree_v$. C\ref{C1} implies that the subgroup $A$
stabilises a vertex $\Lambda$ of $\tree_v$, and hence condition S\ref{S1} is
satisfied. 

Note that $A$ is finitely generated. Condition C\ref{C4} yields that the
action of $\pi_1(M)$ on $\tree_v$ is non--trivial, and Corollary 2 in Section
I.6.5 of \cite{serre} implies
that either $t$ is loxodromic with respect to the action on $\tree_v$ or
there is $a \in A$ such that $ta$ or $at$ is loxodromic.
In the first case, we keep $A_1$ and $A_2$
as they are; in the second case, we replace $t$ by $ta$ and $A_2$ by
$a^{-1}A_2a$; and in the third case, we replace $t$ by $at$ and $A_1$ by
$aA_1 a^{-1}$. Thus, $t$ satisfies condition S\ref{S4}.

Since $A$ stabilises $\Lambda$, $t^{-1}At$ stabilises $t^{-1} \Lambda$, and
since $t$ acts as a loxodromic, we have $t^{-1} \Lambda \ne \Lambda$. In
particular, $A_2$ fixes these two distinct vertices, and hence the path
$[\Lambda, t^{-1} \Lambda]$ pointwise, which implies that it is contained in
an edge stabiliser. Thus, condition S\ref{S2} is satisfied, and the lemma
is proven in the case where $S$ is connected, essential and non--separating,
since condition S\ref{S3} does not apply. 

The proof for the separating case is similar, and will therefore be omitted.
\end{proof}

The conditions are not sufficient when $S$ has more than one
component, since condition C\ref{C4} does not rule out the possibility that
the limiting character is finite on all components of $M-S'$ for a proper
subsurface $S'$ of $S$.


\subsection{The Kinoshita--Terasaka knot}

Let $M$ and $M^\tau$ denote the complements of the Kinoshita--Terasaka knot
and its Conway mutant respectively, and $S$ the corresponding Conway sphere.
In \cite{tillus}, the 
$S$--reducible--non--abelian representations in $\R (M)$ and $\R (M^\tau)$
are computed up to conjugacy, and a comparison thereof leads to the
conclusion that a closed essential surface in $M$ is
associated to an ideal point of $\X (M)$. 

Lemma \ref{associating lem} together with the calculations in \cite{tillus}
implies that the Conway sphere as well as any 
surface obtained by joining boundary components of the sphere with
annuli is a surface associated to the ideal points of $\X (M)$ at which the
holes in the eigenvalue variety occur. Two detected genus two surfaces
and their involutions are shown in Figures 3(a) and 3(b) in \cite{tillus}. 


\subsection{The figure eight knot}

The complement $M$ of the figure eight knot $\k$ in $S^3$ is a once--punctured
torus bundle with fibre a Seifert surface of the knot. Mutation along this
surface results in the so--called sister manifold. The mutation is detected by
the first homology group, but also by the $\SL$--Dehn surgery components. We
verify that the $\PSL$--Dehn surgery components are birationally equivalent,
and use the mutation map to show that the fibres in both manifolds are
detected by ideal points. This method may have non--trivial applications.
\begin{figure}[h]
\psfrag{a}{{\small $a$}}
\psfrag{b}{{\small $b$}}
\psfrag{t}{{\small $t$}}
\psfrag{T}{{\small $\tau$}}
\begin{center}
   \includegraphics[height=4cm]{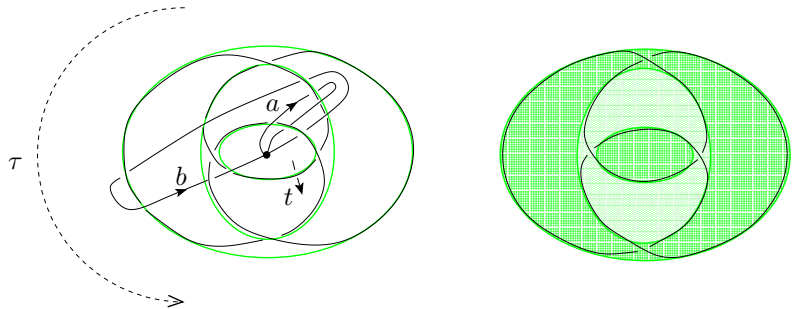}    
  \end{center}
  \caption{Mutation along the Seifert surface}
  \label{fig:fig8}
\end{figure}

A Seifert surface $T_1$ is shown in Figure \ref{fig:fig8}.
A base point and generators are chosen such that
$\tau (a) = a^{-1}$, $\tau (b) = ab^{-1}a^{-1}$, and we compute the
presentation $\G = \langle t,a,b \mid t^{-1}at = aba, t^{-1}bt = ba \rangle$
for $\pi_1(M)$.
The action of $t$ corresponds to the isomorphism $\Phi$ induced by the
monodromy of the fibre bundle.
The isomorphism for the mutative manifold $M^\tau$
is
$\Phi^\tau (a) := \tau(\Phi(a)) = b^{-1}a^{-2}$ and
$\Phi^\tau (b) := \tau(\Phi(b)) = ab^{-1}a^{-2}$,
which yields a presentation $\G^\tau$ for $\pi_1(M^{\tau})$.
Both presentations can be simplified to:
\begin{equation*}
\label{fig8:fund1}
\G = \langle t,a \mid t^{-1}a^{-1}t^{-1}a t a^{-2} t a =1\rangle
\quad\text{and}\quad
\G^\tau = \langle t,a \mid t^{-1}a t a^2 t a t^{-1} a=1\rangle.
\end{equation*}
Note that $H_1(M) \cong \Z$ and $H_1(M^\tau)\cong \Z_5 \oplus \Z$.
Let $x=\tr\rho(t)$ and $y=\tr\rho(a)$. A computation reveals
$\X^{r} (M) = \{ (x,y) \in \C^2 \mid 0 = (2-y)\}$
and  
$\X^{r} (M^\tau) = \{ (x,y) \in \C^2 \mid 0 = (2-y)(1-y-y^2)\}$.

It turns out that the character varieties have only one component
containing the character of an irreducible representation:
\begin{align*}
  \X_0(M) &= \{(x,y)\in \C^2 \mid 0 = 1-y-y^2 +(-1+y)x^2 \},\\
  \X_0(M^\tau) &= \{  (x,y)\in\C^2 \mid 0 = 1+ (-1+y)x^2 \}.
\end{align*}
The curve $\X_{0}(M)$ has no singularities and no singularities at
infinity. Its smooth projective completion is therefore a torus.
The curve $\X_{0}(M^{\tau})$ is rational, and a smooth projective model is
hence a sphere. 

Each $\PSL$--representation lifts to $\SL$ for each example, and 
the quotient map is given by 
$\qe(x,y) = (x^2, y)$. Thus:
\begin{align*}
  \PX (M) &= \{ (X,y) \in \C^2 \mid 0 = (2-y)(1-y-y^2 +(-1+y)X)       \} \\
  \PX (M^\tau) &= \{ (X,y) \in \C^2 \mid 0 = (2-y)(1-y-y^2)(1+(-1+y)X)\}  
\end{align*}
The rational maps between the Dehn surgery components induced by mutation
show that $\PX_0(M)$ and $\PX_0(M^\tau)$ are in fact homeomorphic:
\begin{align*}
\pmap &: \PX_0(M) \to \PX_0(M^\tau ) &&\qquad (X,y) \to 
\big(\frac{1}{1-y},y\big),\\
\pmap^{-1} &: \PX_0(M^\tau ) \to \PX_0(M) &&\qquad (X,y) \to 
\big(\frac{1-y-y^2}{1-y},y\big).
\end{align*}

The surfaces detected by the Dehn surgery components do not include the
fibre, but one can recover curves of reducible representations as
follows. There 
are only three points on each of the projective Dehn surgery components on
which $\pmap$ and $\pmap^{-1}$ are not defined a
priori, and they correspond to the intersection with
$\{(2-y)(1-y-y^2)=0\}$. 
The corresponding representations of $M$ and $M^{\tau}$ are 
$T_1$--abelian and
satisfy $\tau (\rho (\gamma )) = \rho (\gamma)^{-1}$. For each we
can find a
1--parameter family of elements in $\PSL$ which realise the action of $\tau$.
Consider the following lift to $\SL$ of an irreducible
$\PSL$--representation of $M$: 
\begin{equation*}
  \rho (t) =  \begin{pmatrix} i & 1 \\ 0 & i^{-1} \end{pmatrix}
  \quad \text{and} \quad
  \rho (a) = \begin{pmatrix} u & 0 \\ i(u^{-1}-u) & u^{-1}
\end{pmatrix}
\end{equation*}
subject to $0=1+u+u^2+u^3+u^4$. These are dihedral representations, and they
are abelian on the fibre. Elements realising the involution are: 
\begin{equation*}
  H(z) = \begin{pmatrix} iz & z \\ z-z^{-1} & i^{-1}z \end{pmatrix}
  \quad \text{for any $z \in \C-\{0\}$},
\end{equation*}
and we obtain the following representations $\rho_z \in \R (M^\tau)$:
\begin{equation*}
  \rho_z (t) =H(z)\rho (t) = \begin{pmatrix} -z & 0 \\ i(z-z^{-1}) & -z^{-1}
  \end{pmatrix} 
  \thinspace \text{and} \thinspace\thinspace
  \rho_z (a) =\rho (a).
\end{equation*}
These representations are abelian. The construction
yields a map $\C-\{0\} \to \PX (M^\tau)$, which is non--constant
since $(\tr\rho_z (t))^2 = (z + z^{-1})^2$, and the image is a
curve in $\PX (M^\tau)$. At an ideal point of this curve, the conditions of
Lemma \ref{associating lem} are satisfied with respect to the fibre in
$M^{\tau}$. One can do a similar construction for the other points in
$\PX^i (M) \cap \PF (M)$. 

Using characters in
$\PX^i (M^\tau) \cap \PF (M^\tau)$, one only obtains a curve in $\PX (M)$
for the point corresponding to the intersection with $\{ y=2 \}$, whilst the
points in the intersection with $\{ 1=y+y^2\}$ yield a map
$\C \to \PX (M)$ which is constant.


\Addresses\recd

\begin{thebibliography}

\bibitem{bozh} S. Boyer, X. Zhang: \emph{On Culler--Shalen seminorms
   and Dehn filling}, Ann. Math., 148, 737-801 (1998).

\bibitem{a-poly1} D. Cooper, D.D. Long: \emph{Remarks on the A-polynomial of
  a Knot}, J. Knot Theory Ramifications, 5, 609-628 (1996).

\bibitem{cs} M. Culler, P.B.~Shalen: \emph{Varieties of group
  representations and splittings of 3-manifolds}, Ann. Math., 117, 109-146
      (1983).

\bibitem{g} W. Goldman:
  \emph{Topological components of spaces of representations},
      Invent. Math., 93, 557-607 (1988).

\bibitem{r} D. Ruberman: \emph{Mutation and volumes of knots in $S^3$},
          Invent. Math., 90, 189-215 (1987).

\bibitem{serre} J.--P. Serre: \emph{Trees}, Springer, Berlin,
    1980.

\bibitem{sh1} P.B. Shalen: \emph{Representations of 3--manifold groups},
Handbook of geometric topology, 955-1044, Amsterdam, 2002.


\bibitem{t} W.P. Thurston: \emph{The geometry and topology of
3--manifolds}, Princeton Univ. Math. Dept. (1978). 
\url{http://msri.org/publications/books/gt3m/}

\bibitem{tillus} S. Tillmann: \emph{On the Kinoshita--Terasaka mutants
    and generalised Conway mutation}, J. Knot Theory Ramifications, 9,
    557-575 (2000). 

\end{thebibliography}
\end{document}